# COMPACTNESS IN THE $\overline{\partial}$-NEUMANN PROBLEM

SIQI FU AND EMIL J. STRAUBE

## Contents

## 1. Introduction

Let $\Omega$ be a bounded pseudoconvex domain in $\mathbb{C}^n$. The $\overline{\partial}$ (or Dolbeault)-complex is the complex

$$(1) \quad L^2(\Omega) \xrightarrow{\overline{\partial}} L^2_{(0,1)}(\Omega) \xrightarrow{\overline{\partial}} L^2_{(0,2)}(\Omega) \xrightarrow{\overline{\partial}} \cdots \xrightarrow{\overline{\partial}} L^2_{(0,n)}(\Omega) \xrightarrow{\overline{\partial}} 0,$$

where $L^2_{(0,q)}(\Omega)$ denotes the space of $(0,q)$-forms on $\Omega$ with coefficients in $L^2(\Omega)$. The norm (with corresponding inner product) is $\|\sum'_J a_J\, d\bar{z}_J\|^2 = \sum'_J \int_\Omega |a_J|^2 dV$. Here, the prime denotes summation over strictly increasing $q$-tuples $J$, and $d\bar{z}_J = d\bar{z}_{j_1} \wedge d\bar{z}_{j_2} \wedge \cdots \wedge d\bar{z}_{j_q}$. The $\overline{\partial}$-operator acts via

$$(2) \quad \overline{\partial}\left(\sum_J{}' a_J\, d\bar{z}_J\right) = \sum_{j=1}^n \sum_J{}' \frac{\partial a_J}{\partial \bar{z}_j} d\bar{z}_j \wedge d\bar{z}_J.$$

The derivatives are taken in the distribution sense, and the domain of $\overline{\partial}$ consists of those $(0,q)$-forms where the right hand side is in $L^2_{(0,q+1)}(\Omega)$. $\overline{\partial}$ is then a densely defined closed operator, and so has an adjoint mapping $L^2_{(0,q+1)}(\Omega)$ into $L^2_{(0,q)}(\Omega)$. We denote this adjoint by $\overline{\partial}^*$. The complex Laplacian is the operator $\square = \overline{\partial}\overline{\partial}^* + \overline{\partial}^*\overline{\partial}$, acting as an (unbounded) operator on $L^2_{(0,q)}(\Omega)$.

The basic $L^2$-existence theorem for the $\overline{\partial}$-Neumann problem (the problem of inverting $\square$) goes back to Hörmander[H65]: for each $q$, $1 \leq q \leq n$, $\square = \overline{\partial}\overline{\partial}^* + \overline{\partial}^*\overline{\partial}$ is an unbounded self-adjoint operator on $L^2_{(0,q)}(\Omega)$ which is surjective and (consequently) has a bounded inverse. This inverse is the $\overline{\partial}$-Neumann operator $N_q$.

1991 *Mathematics Subject Classification.* 32F20, 35N15, 32F15.
This work is supported in part by NSF grants DMS-9632028 and DMS-9801539.

More precisely, $N_q$ satisfies the estimate

$$\text{(3)} \qquad \|N_q u\| \leq \left(\frac{D^2 e}{q}\right) \|u\|, \quad u \in L^2_{(0,q)}(\Omega),$$

where $D$ is the diameter of the domain $\Omega$. If $\alpha$ is a $\overline{\partial}$-closed $(0, q+1)$-form, then $\overline{\partial}^* N_{q+1} \alpha$ provides the solution orthogonal to the kernel of $\overline{\partial}$ (*i.e.*, the solution of minimal norm) to $\overline{\partial} u = \alpha$. The $\overline{\partial}$-Neumann operator is also closely related to the Bergman projection $P_q$, the orthogonal projection from $L^2_{(0,q)}$ onto its subspace consisting of $\overline{\partial}$-closed forms, via Kohn's formula $P_q = I - \overline{\partial}^* N_{q+1} \overline{\partial}$. We refer the reader to [FK72, Kr92, BS99, CS99] for further details and references on the $L^2$-theory of the $\overline{\partial}$-Neumann problem.

The question we consider in this article is that of compactness of $N_q$. This question is of interest for various reasons. We indicate some consequences of compactness in Section 2 below. In addition, from the point of view of finding necessary and sufficient conditions in terms of properties of the boundary, compactness appears to be more tractable than global regularity.

The compactness condition can be reformulated in several useful ways.

**Lemma 1.** *Let $\Omega$ be a bounded pseudoconvex domain, $1 \leq q \leq n$. Then the following are equivalent:*

1. *The $\overline{\partial}$-Neumann operator $N_q$ is compact from $L^2_{(0,q)}(\Omega)$ to itself.*
2. *The embedding of the space $\operatorname{Dom} \overline{\partial} \cap \operatorname{Dom} \overline{\partial}^*$, provided with the graph norm $u \to \|\overline{\partial} u\| + \|\overline{\partial}^* u\|$, into $L^2_{(0,q)}(\Omega)$ is compact.*
3. *For every $\epsilon > 0$ there exists a constant $C_\epsilon > 0$ such that*

$$\|u\|^2 \leq \epsilon \left(\|\overline{\partial} u\|^2 + \|\overline{\partial}^* u\|^2\right) + C_\epsilon \|u\|^2_{-1}$$

   *when $u \in \operatorname{Dom} \overline{\partial} \cap \operatorname{Dom} \overline{\partial}^*$.*
4. *The canonical solution operators $\overline{\partial}^* N_q \colon L^2_{(0,q)}(\Omega) \to L^2_{(0,q-1)}(\Omega)$ and $\overline{\partial}^* N_{q+1} \colon L^2_{(0,q+1)}(\Omega) \to L^2_{(0,q)}(\Omega)$ are compact.*

The statement in (3) is called a compactness estimate. In (4), we refer to $\overline{\partial}^* N_q$ and $\overline{\partial}^* N_{q+1}$ as "canonical solution operators", although strictly speaking the solution operators are the restrictions to the kernel of $\overline{\partial}$. However, $\overline{\partial}^* N_q$ and $\overline{\partial}^* N_{q+1}$ are zero on the orthogonal complement of $\ker \overline{\partial}_q$ and $\ker \overline{\partial}_{q+1}$, respectively. In particular, compactness is not affected. Note that saying that the canonical solution operator is compact is the same as saying that there exists *some* compact solution operator, as the projection onto the orthogonal complement of the kernel of $\overline{\partial}$ preserves compactness. The equivalence of (2) and (3) is in [KN65], Lemma 1.1. That (1) is equivalent to (2) and (3) is an easy consequence of the general $L^2$-theory and the fact that $L^2_{(0,q)}(\Omega)$ embeds compactly into $W^{-1}_{(0,q)}(\Omega)$. Finally, the equivalence of (1) and (4) follows from the formula $N_q = (\overline{\partial}^* N_q)^*(\overline{\partial}^* N_q) + (\overline{\partial}^* N_{q+1})(\overline{\partial}^* N_{q+1})^*$ (see [FK72], p.55, [Ra84]).

It is useful to know that compactness is a *local* property. Roughly speaking, the $\overline{\partial}$-Neumann operator $N_q$ on $\Omega$ is compact if and only if every boundary point has a neighborhood $U$ such that the corresponding $\overline{\partial}$-Neumann operator on $U \cap \Omega$ is compact. For simplicity, we will assume that $\Omega$ and $U$ are such that these intersections are domains (*i.e.*, connected); this is not a problem in the applications we have in mind (see section 5).

**Lemma 2.** *Let $\Omega$ be a bounded pseudoconvex domain in $\mathbb{C}^n$, $1 \leq q \leq n$.*

1. *If for every boundary point there exists a pseudoconvex domain $U$ such that $N_q$ on (the domain) $U \cap \Omega$ is compact, then $N_q$ on $\Omega$ is compact.*
2. *If $U$ is smoothly bounded and strictly pseudoconvex and $U \cap \Omega$ is a domain, then if the $\overline{\partial}$-Neumann operator on $(0,q)$-forms on $\Omega$ is compact, so is the corresponding $\overline{\partial}$-Neumann operator on $U \cap \Omega$.*

The proof of (1) is straightforward and results from a partition of unity argument together with the interior elliptic regularity of $\overline{\partial} \oplus \overline{\partial}^*$. The proof of (2) is also standard, but it needs some ideas discussed in Section 3 below. Accordingly, we postpone this proof; it will be given at the end of Section 3.

## 2. Some consequences of compactness

In the case of smoothly bounded domains, a classical theorem of Kohn and Nirenberg [KN65] asserts that compactness of $N_q$ (as an operator from $L^2_{(0,q)}(\Omega)$ to itself) implies global regularity in the sense of preservation of Sobolev spaces. We denote by $W^s(\Omega)$ the standard $L^2$-Sobolev spaces on $\Omega$, and by $W^s_{(0,q)}(\Omega)$ the spaces of $(0,q)$-forms on $\Omega$ with coefficients in $W^s(\Omega)$.

**Theorem 3.** *Let $\Omega$ be a bounded pseudoconvex domain in $\mathbb{C}^n$ with smooth boundary. Let $1 \leq q \leq n$. If $N_q$ is compact on $L^2_{(0,q)}(\Omega)$, then $N_q$ is compact (in particular, continuous) as an operator from $W^s_{(0,q)}(\Omega)$ to itself, for all $s \geq 0$.*

**Remark 1.** In Theorem 3, the implication in the other direction is valid as well: if $N_q$ is compact as an operator from $W^s_{(0,q)}(\Omega)$ to itself for *some* $s \geq 0$, then $N_q$ is compact in $L^2_{(0,q)}(\Omega)$ (and hence in $W^s_{(0,q)}(\Omega)$, for all $s \geq 0$). This follows from an abstract theorem about compact operators over a Banach space which are symmetric with respect to a scalar product. Such operators are necessarily compact with respect to the norm given by the scalar product (see, e.g., [La54], Corollary 2).

The Fredholm theory of Toeplitz operators is an immediate consequence of compactness of the $\overline{\partial}$-Neumann problem [V72, HI97]. In fact, compactness of the canonical solution operators implies that commutators between the Bergman projection and multiplication operators are compact [CD97]:

**Proposition 4.** *Let $\Omega$ be a bounded pseudoconvex domain in $\mathbb{C}^n$. Assume that for some $q$, $0 \leq q \leq n-1$, the canonical solution operator $\overline{\partial}^* N_{q+1}$ is compact. Let $M$ be a function that has bounded first order partial derivatives on $\Omega$. Then the commutators $[P_q, M]$ between the Bergman projection $P_q$ and the multiplication operator by $M$ is compact on $L^2_{(0,q)}(\Omega)$.*

A short proof is as follows. For $f \in L^2_{(0,q)}(\Omega)$, let $f_1 := P_q f$, $f_2 := (I - P_q)f$. Since both projections are continuous, it suffices to see that $[P_q, M]$ is compact on both $\ker \overline{\partial}$ and $(\ker \overline{\partial})^\perp$. We have

$$
\begin{aligned}
[P_q, M]f_1 &= P_q M f_1 - M P_q f_1 = M f_1 - \overline{\partial}^* N_{q+1} \overline{\partial} M f_1 - M f_1 \\
&= -\overline{\partial}^* N_{q+1} \overline{\partial}(M f_1) = -\overline{\partial}^* N_{q+1}(\overline{\partial} M \wedge f_1).
\end{aligned}
\tag{4}
$$

(Compare [Sa91], proof of Lemma 3 for a similar argument.) Since $\overline{\partial}M$ is bounded in $\Omega$ (so acts as an $L^2$-bounded multiplier) and since $\overline{\partial}^* N_{q+1}$ is compact, $[P_q, M]$ is compact on $\ker \overline{\partial}$. On the other hand,

$$
\begin{aligned}
[P_q,\ M]f_2 &= P_q M f_2 - M P_q f_2 = P_q M f_2 \\
&= P_q M \overline{\partial}^* N_{q+1} \overline{\partial} f_2 = P_q [M,\ \overline{\partial}^*](N_{q+1}\overline{\partial})f_2.
\end{aligned}
\tag{5}
$$

Here we have used the fact that $P_q \overline{\partial}^* M N_{q+1} \overline{\partial} f_2 = 0$, since the range of $\overline{\partial}^*$ is orthogonal to $\overline{\partial}$-closed forms (note that $N_{q+1}$ maps into the domain of $\overline{\partial}^*$, which is preserved by multiplication with $M$). Since $[M, \overline{\partial}^*]$ acts as a zero order (bounded) operator on $N_{q+1}\overline{\partial} f_2$, and $N_{q+1}\overline{\partial}$ is compact (since $N_{q+1}\overline{\partial} = (\overline{\partial}^* N_{q+1})^*$), we are done.

**Remark 2.** In light of Lemma 1, (4), and Proposition 4, compactness of either $N_q$ or $N_{q+1}$ implies the compactness of $[P_q, M]$. It would be of interest to see to what extent properties of commutators between the Bergman projections and multiplication operators can be used to actually characterize compactness properties of the $\overline{\partial}$-Neumann problem. Compactness of the commutators $[P_q,\ M]$ does imply compactness of the canonical solution operator $\overline{\partial}^* N_{q+1}$ *restricted* to forms with *holomorphic* coefficients. This is essentially Lemma 2 in [Sa91]. Indeed, suppose that $u \in L^2_{(0,q+1)}(\Omega)$ has holomorphic coefficients. Write $u = \sum_j \sum'_K u_{jK} d\bar{z}_j \wedge d\bar{z}_K$, $1 \leq j \leq n$, $|K| = q$, then by (4) above, $-\sum_j [P_q,\ \bar{z}_j](\sum'_K u_{jK} d\bar{z}_K) = \sum_j \overline{\partial}^* N_{q+1}(\sum'_K u_{jK} d\bar{z}_j \wedge d\bar{z}_K) = \overline{\partial}^* N_{q+1} u$. Therefore, $\overline{\partial}^* N_{q+1}$ is compact on forms with holomorphic coefficients if the commutators $[P_q, \bar{z}_j]$, $1 \leq j \leq n$, are. On *convex* domains, compactness of this restriction *is* sufficient to give compactness of $N_{q+1}$, see [FS98], Remark (2), Section 5.

Whether or not the $\overline{\partial}$-Neumann problem is compact has further ramifications for the theory of the Toeplitz $C^*$-algebras naturally associated to a domain, see [Sa91, Sa95].

## 3. Sufficient conditions for compactness

We start with an inequality due to Catlin that is central to the subject (see [Ca84a], p. 45; [Ca87], Theorem 2.1; see also [BS99], Section 2, for a somewhat different approach to this type of estimate). Assume for the moment that the boundary of $\Omega$ is sufficiently smooth, say $C^2$, to allow integration by parts. Let $\lambda \in C^2(\overline{\Omega})$ be normalized so that $0 \leq \lambda \leq 1$. Then

$$
\sum_K{}' \sum_{j,k=1}^n \int_\Omega \frac{\partial^2 \lambda}{\partial z_j \partial \bar{z}_k} u_{jK} \overline{u}_{kK} dV \leq e\left(\|\overline{\partial} u\|^2 + \|\overline{\partial}^* u\|^2\right)
\tag{6}
$$

for all $u \in \operatorname{Dom}\overline{\partial} \cap \operatorname{Dom}\overline{\partial}^* \subset L^2_{(0,q)}(\Omega)$. Heuristically, in view of (3) in Lemma 1 above (the compactness estimate), (6) shows that functions $\lambda$ with "large" Hessians are desirable when compactness is the goal. To make this precise, it is convenient to recall the following notions (see [Ca87]).

Let $\Lambda_z^{(0,q)}$ be the space of $(0,q)$-forms at $z$ equipped with the standard Hermitian metric $|\sum'_J u_J d\bar{z}_J|^2 = \sum'_J |u_J|^2$. For a $C^2$-smooth function $\lambda(z)$ in a neighborhood

of $z$, let

$$H_q(\lambda)(z,u) = {\sum_{|K|=q-1}}' \sum_{j,k=1}^{n} \frac{\partial^2 \lambda(z)}{\partial z_j \partial \bar{z}_k} u_{jK} \bar{u}_{kK}, \qquad u \in \Lambda_z^{(0,q)}.$$

Note that $H_1(\lambda)$ is the usual complex Hessian.

**Lemma 5.** *Let $M \geq 0$ be a constant. The following statements are equivalent.*

1. $H_q(\lambda)(z,u) \geq M|u|^2$ for all $u \in \Lambda_z^{(0,q)}$.
2. *For any orthonormal subset $\{t^j, \; 1 \leq j \leq q\}$ of $\mathbb{C}^n$,*

$$\sum_{j=1}^{q} H_1(\lambda)(z,t^j) \geq M.$$

3. *The sum of any $q$ eigenvalues of the Hermitian matrix $(\partial^2 \lambda / \partial z_j \partial \bar{z}_k)$ is greater than or equal to $M$.*

The equivalence of (1) and (2) is in [Ca87], pp. 189-190. The equivalence of (1) and (3) is implicit in [H65], pp. 137. This lemma is proved by first establishing it in the case when the Hermitian matrix $(\partial^2 \lambda / \partial z_j \partial \bar{z}_k)$ is diagonal, and then noting that statements (1)-(3) are invariant under a unitary transformation.

In light of Lemma 5, it is useful to have the following definition: For a bounded pseudoconvex domain $\Omega$, we say that $b\Omega$ satisfies *property ($P_q$)* (or just ($P_q$)) if for every positive number $M$, there exists a neighborhood $U$ of $b\Omega$ and a $C^2$-smooth function $\lambda$ on $U \cap \Omega$, $0 \leq \lambda \leq 1$, such that for all $z$ in $U \cap \Omega$, the sum of any $q$ (equivalently: the smallest $q$) eigenvalues of the Hermitian form $(\partial \lambda / \partial z_j \partial \bar{z}_k)_{j,k=1}^{n}$ is at least $M$. This definition is from [Ca84a] (for $q = 1$). Note that $(P_1) \Rightarrow (P_2) \Rightarrow \cdots \Rightarrow (P_n)$. With this definition, we can make our heuristic statement from above precise.

**Theorem 6.** *Let $\Omega$ be a bounded pseudoconvex domain in $\mathbb{C}^n$, $1 \leq q \leq n$. If the boundary of $\Omega$ satisfies property ($P_q$), then the $\bar{\partial}$-Neumann operator $N_q$ is compact.*

For sufficiently smooth domains, Theorem 6 is in [Ca84a], proof of Theorem 1 (see also [Ca87]). The boundary regularity requirement was considerably weakened in [HI97], Theorem 1. The second author showed in [St97], Corollary 3, that no boundary regularity at all is needed.

**Remark 3.** The differentiability requirement on $\lambda$ can be relaxed. For example, assume only that $\lambda$ is continuous with its complex Hessian, interpreted as a *current*, bounded below by $M dd^c |z|^2$. Then $N_q$ is compact, $1 \leq q \leq n$. This is also contained in [St97].

The proof of Theorem 6 results from the following considerations. Assume first that $b\Omega$ is smooth, so that we can apply (6). We note that by Lemma 5, the condition on $\lambda$ (the function in the definition of property ($P_q$)) implies that

(7) $$M|u(z)|^2 \leq H_q(\lambda)(z,u)$$

for all $(0,q)$-forms $u$ and $z \in U \cap \Omega$. Choosing $M \approx 1/\epsilon$, (7) combines with (6) to produce the desired compactness estimate, modulo terms that are compactly supported ((7) holds only near the boundary). These latter terms, however, are easily handled by the interior elliptic regularity of $\bar{\partial} \oplus \bar{\partial}^*$ (see [Ca84a]).

When no boundary regularity is assumed, one would like to obtain (6) by first integrating over approximating subdomains, and then passing to the limit. The reason that this is not immediate is that when forms are restricted to subdomains, they are in general no longer in the domain of $\overline{\partial}^*$ on the subdomain. This difficulty can be overcome by a regularization procedure introduced in [St97] that exploits the $\overline{\partial}$-Neumann operators of the subdomains. We refer the reader to [St97], proof of Corollary 3, for details.

The simplest examples of domains satisfying $(P_1)$ (hence $(P_q)$, $1 < q$) are strictly pseudoconvex domains: it suffices to consider a strictly plurisubharmonic defining function. More generally, domains of finite type satisfy $(P_1)$. This is far from obvious, however. It is a consequence of the analysis of finite type points in [Ca84b] and [Da82], see [Ca84a]. But property $(P_1)$ is considerably more general. It is not hard to see that if a domain is strictly pseudoconvex except at finitely many points, then it satisfies $(P_1)$. In fact, if the infinite type points of the boundary have 2-dimensional Hausdorff measure zero, then the domain satisfies $(P_1)$ ([Bo88]; [Si87], Remark on p. 310). On the other hand, Sibony found examples of domains with $(P_1)$ where the set of boundary points of infinite type has positive measure ([Si87], p. 310).

Sibony undertook a systematic study of property $(P_1)$ in the context of general compact sets in $\mathbb{C}^n$. One of the main tools is Choquet theory applied to the cone of plurisubharmonic functions. One can carry out an analogous study of property $(P_q)$ by considering a cone of functions that reflects the condition on the Hessian used in Lemma 5, as in [FS98]. For an open set $U \subset \mathbb{C}^n$, denote by $P_q(U)$ the set of continuous functions $\lambda$ on $U$ such that for any $z \in U$ and orthonormal set of vectors $\{t_1, \cdots, t_q\}$ in $\mathbb{C}^n$, the function

$$\zeta = (\zeta_1, \cdots, \zeta_q) \in \mathbb{C}^q \mapsto \lambda(z + \zeta_1 t_1 + \ldots + \zeta_q t_q)$$

is subharmonic on $\{\zeta \in \mathbb{C}^q;\ z + \zeta_1 t_1 + \ldots + \zeta_q t_q \in U\}$. That is, $P_q(U)$ consists of the continuous functions on $U$ that are subharmonic on each $q$-dimensional complex affine subspace. In particular, $P_1(U)$ is the set of all continuous plurisubharmonic functions and $P_n(U)$ is the set of all continuous subharmonic functions. $P_q(U)$ is a convex cone in $C(U)$ that is closed under taking the pointwise maximum of finitely many of its elements. Each function in $P_q(U)$ is a locally uniform limit of $C^\infty$-smooth elements in $P_q$ of slightly smaller open sets: this follows from the usual mollifier argument. It is easy to check that $-\sum_{j=1}^{q-1} |z_j|^2 + (q-1) \sum_{j=q}^n |z_j|^2 \in P_q(\mathbb{C}^n)$.

For a compact subset $X$ of $\mathbb{C}^n$, let $P_q(X)$ be the uniform closure in $C(X)$ of functions that are in $P_q(U)$ for some neighborhood $U$ of $X$. A probability measure $\mu$ on $X$ is said to be a $P_q$-measure for $z \in X$ if

$$\lambda(z) \le \int_X \lambda\, d\mu, \qquad \text{for all}\ \ \lambda \in P_q(X).$$

Let $J_q(X)$ be the Choquet boundary of $P_q(X)$, i.e., the subset of $X$ consisting of points $z \in X$ such that the point mass is the only $P_q$-measure for $z$. We refer the reader to [G78], Chapter 1 for the elements of Choquet theory. In particular, characterizations of the Choquet boundary are given in Theorem 1.13 in [G78]. We are interested in the situation where the Choquet boundary is all of $X$.

**Proposition 7.** *Let $X$ be a compact subset of $\mathbb{C}^n$. Then the following are equivalent:*

1. $J_q(X) = X$.
2. $C(X) = P_q(X)$.
3. *For any $z \in X$, there exists an $r > 0$ such that if $X_r = X \cap \overline{B(z,r)}$ then $C(X_r) = P_q(X_r)$.*
4. *For any $M > 0$, there exists a function $\lambda$ smooth in a neighborhood of $X$ such that $0 < \lambda(z) < 1$ and $H_q(\lambda)(z,u) \geq M|u|^2$ for $z \in X$ and $u \in \Lambda_z^{(0,q)}$.*

Proposition 7 is from [Si87] (Proposition 1.3 and 1.4, compare also [Si89]) when $q = 1$, but the arguments cover the case of general $q$ as well (in fact, many of the arguments hold in the general context of [G78], Chapter 1). When $q = 1$, Sibony calls a compact set that has the (equivalent) properties in Proposition 7 $B$-regular. Accordingly, we shall use the term $B_q$-*regular* when $q \geq 1$. The significance of Proposition 7 in our context stems from (4): when $X = b\Omega$, where $\Omega$ is a bounded domain, then (4) is essentially property $(P_q)$ for $b\Omega$. Actually, property $(P_q)$ requires the function $\lambda$ to exist only in a neighborhood of $b\Omega$ intersected with $\overline{\Omega}$, rather than in a full neighborhood of $b\Omega$. However, it is easy to see that on domains with relatively minimally regular boundary (for instance, when the boundary is locally a graph), the two notions coincide.

**Remark 4.** For boundaries of bounded domains which admit a bounded plurisubharmonic exhaustion function (hyperconvex domains), having the (equivalent) properties in Proposition 7 can also be characterized in terms of various properties related to the potential theory of the cone $P_q(\Omega)$, for example the existence of peak functions in $P_q(\Omega) \cap C(\overline{\Omega})$. For this, see [Si87], Theorem 2.1, and [Si89], Theorem 2.3, where the case $q = 1$ is treated. However, Sibony's arguments carry over to the case of general $q$.

The easiest way that $B$-regularity can fail is for the boundary to contain an analytic disc: in view of the maximum principle, $C(b\Omega) = P_q(b\Omega)$ cannot hold. A more direct argument is to pull back to the unit disc and to observe that subharmonic functions on the unit disc, with values between 0 and 1, cannot have arbitrarily large Hessians. Similarly, $q$-dimensional varieties in the boundary are incompatible with $B_q$-regularity. However, the absence of varieties is not sufficient for $B_q$-regularity in general: Sibony has given examples of complete smooth Hartogs domains in $\mathbb{C}^2$ whose boundaries contain no discs, yet are not $B$-regular. ([Si87], p. 310). It turns out that on locally convexifiable domains, $B_q$-regularity (hence property $(P_q)$) *is* equivalent to the absence of $q$-dimensional varieties from the boundary; we will discuss this in section 5 below.

$B_q$-regularity is preserved under countable unions.

**Proposition 8.** *Let $X_k$, $k = 1, 2, 3, \ldots$, be compact subsets of $\mathbb{C}^n$, all of which are $B_q$-regular. If $X = \cup_{k=1}^{\infty} X_k$ is compact, then $X$ is $B_q$-regular.*

For $q = 1$, this is Proposition 1.9 in [Si87]. Sibony's proof works essentially verbatim to give the result for all $q$. (This includes the supporting Propositions 1.4 and 1.6, and Lemma 1.8.)

We conclude this section by completing the proof of Lemma 2 from Section 1. Assume that there is a compactness estimate for $(0, q)$-forms on $\Omega$, and let $U$ be strictly pseudoconvex (and $U \cap \Omega$ a domain). We show how to obtain a compactness estimate on $U \cap \Omega$. Let $u \in \operatorname{Dom}\overline{\partial} \cap \operatorname{Dom}\overline{\partial}^* \subset L^2_{(0,q)}(U \cap \Omega)$. Fix $\epsilon > 0$. The boundary of $U$ is strictly pseudoconvex, so satisfies property $(P_1)$, hence property $(P_q)$ for all $q$. Thus for a smooth cutoff function $\phi_\epsilon$ that is identically equal to 1 on $bU$ and supported in a sufficiently small neighborhood of $bU$, we get by the arguments in the proof of Theorem 6:

$$\begin{aligned} \|u\|^2_{U\cap\Omega} &\lesssim \|\phi_\epsilon u\|^2_{U\cap\Omega} + \|(1-\phi_\epsilon)u\|^2_{U\cap\Omega} \\ &\lesssim \epsilon\left(\|\overline{\partial}(\phi_\epsilon u)\|^2_{U\cap\Omega} + \|\overline{\partial}^*(\phi_\epsilon u)\|^2_{U\cap\Omega}\right) + \|(1-\phi_\epsilon)u\|^2_{U\cap\Omega} \\ &\lesssim \epsilon\left(\|\overline{\partial}u\|^2_{U\cap\Omega} + \|\overline{\partial}^* u\|^2_{U\cap\Omega}\right) + C_\epsilon\|\chi_\epsilon u\|^2_{U\cap\Omega} + \|(1-\phi_\epsilon)u\|^2_{U\cap\Omega}. \end{aligned}$$

Here, $\chi_\epsilon$ is a smooth cutoff function identically 1 on the support of $\bigtriangledown\phi_\epsilon$ whose support is compact in $U$, and $C_\epsilon := \max|\bigtriangledown\phi_\epsilon|$. Now view the forms $\chi_\epsilon u$ and $(1-\phi_\epsilon)u$ as forms on $\Omega$; they are in the domain of $\overline{\partial}^*$ on $\Omega$. Applying the compactness estimate on $\Omega$, with an $\epsilon'$ *sufficiently small*, yields the desired compactness estimate on $U \cap \Omega$ (terms involving $\|u\|^2$ can be absorbed if $\epsilon'$ is chosen small enough). This completes the proof of Lemma 2.

## 4. Discs in the boundary vs. compactness

In this section, we restrict our attention to compactness of the $\overline{\partial}$-Neumann operator on $(0,1)$-forms. Since a disc in the boundary is a blatant violation of $B$-regularity, and given that discs are known to prevent hypoellipticity ([Ca81], [DP81]), it is natural to ask whether there is a corresponding failure of compactness.

Here is a simple example (taken from [Li85], [Kr88]). Let $\Omega$ be the bidisc in $\mathbb{C}^2$, with the edge rounded. Denote by $D$ the unit disc in $\mathbb{C}$, and let $f \in L^2(D) \cap \ker\overline{\partial}$. Using subscripts to denote $L^2$-norms on domains, we have that $\|\bar{z}_2 f(z_1)\|_\Omega \approx \|f\|_D \approx \|f\|_\Omega$. Also, $\overline{\partial}(\bar{z}_2 f(z_1)) = f(z_1)d\bar{z}_2$, and since $\bar{z}_2 f(z_1)$ is orthogonal in $L^2(\Omega)$ to the holomorphic functions, it is the canonical solution, i.e., $\bar{z}_2 f(z_1) = \overline{\partial}^* N_1(f(z_1)d\bar{z}_2)$. Since the norms compare, compactness of $\overline{\partial}^* N_1$ (hence of $N_1$ by (4) in Lemma 1) would imply compactness of the unit ball in $L^2(D) \cap \ker\overline{\partial}$. Consequently, $N_1$ is not compact.

On the other hand, here is an example of a pseudoconvex Reinhardt domain in $\mathbb{C}^2$ with a disc in the boundary, but with compact $\overline{\partial}$-Neumann operator $N_1$. The domain is incomplete and non-smooth. Let $\Omega = \{(z_1, z_2) \in \mathbb{C}^2;\ |z_1|^2 + |z_2|^2 < 1,\ 0 < |z_1| < 1\}$, i.e., $\Omega$ is the unit ball in $\mathbb{C}^2$ minus the variety $\{z_1 = 0\}$. The point is that the $L^2$-theory does not detect the deletion of this variety, and so is the same as on the unit ball. More precisely, the natural isometry $L^2_{(0,q)}(\Omega) \hookrightarrow L^2_{(0,q)}(B)$ ($B$ denotes the unit ball) commutes with $\overline{\partial}_q$, $q = 0, 1, 2$. This can be checked by an argument completely analogous to that in [Be82], p. 687. But then the $\overline{\partial}^*$ operators have similar commutation properties, hence so do the $\square$ operators and their inverses, the $\overline{\partial}$-Neumann operators. In particular, $N_1$ on $\Omega$ inherits compactness from $N_1$ on $B$.

When some boundary regularity is assumed, the phenomenon in the above example cannot occur on domains in $\mathbb{C}^2$.

**Proposition 9.** *Let $\Omega$ be a bounded pseudoconvex domain in $\mathbb{C}^2$ with Lipschitz boundary. If the boundary of $\Omega$ contains an analytic disc, then the $\overline{\partial}$-Neumann operator $N_1$ on $\Omega$ is not compact.*

As of this writing, it is not known whether Proposition 9 holds in higher dimensions. Proposition 9 has been part of the folklore for many years; it is usually attributed to David Catlin; Michael Christ has also found a proof [Ch98].

A proof of Proposition 9 results from adapting the ideas in [Ca81] and [DP81]. We have also used these ideas in [FS98], Section 4. By Lemma 1, (4), it suffices to show that the canonical solution operator $\overline{\partial}^* N_1$ is not compact. We first assume that the disc in the boundary is an affine disc, say the disc $\overline{D}_{2r_0} \times \{0\}$, where $D_r := \{z \in \mathbb{C}; \ |z| < r\}$. After a suitable complex linear change of coordinates, $\Omega$ is defined near the origin by $y_2 < \rho(z_1, x_2)$ ($z_j = x_j + iy_j$, $j = 1, 2$), where $\rho$ is a Lipschitz function. Since $\rho(z_1, 0) = 0$ for $z_1 \in D_{2r_0}$, $|\rho(z_1, x_2)| = |\rho(z_1, x_2) - \rho(z_1, 0)| \leq C|x_2|$ when $|x_2|$ is sufficiently small. Therefore, there exist circular wedges $W_0$ and $W_1$ of the same radius in the $z_2$-plane, symmetric about the $y_2$-axis, such that $D_{r_0} \times W_0 \subset \Omega$, and $\Omega \cap (D_{r_0} \times D_{r_1}) \subset D_{r_0} \times W_1$, for a suitable $r_1 > 0$ (shrink $r_0$ if necessary). We first observe that for any $r_3 > 0$ smaller than the (common) radius of $W_0$ and $W_1$, there is a sequence of holomorphic functions in $W_1$, bounded in $L^2(W_1)$, so that no subsequence converges in $L^2(W_0 \cap D_{r_3})$; i.e., the restriction operator from the Bergman space on $W_1$ to the Bergman space on $W_0 \cap D_{r_3}$ is not compact. This is easily seen by considering the sequence

$$f_j(z_2) := \sqrt{\frac{2 - 2a_j}{\alpha}} R^{a_j - 1} z_2^{-a_j}$$

(defined via a branch cut along the positive imaginary axis), where $a_j \nearrow 1$. Here, $R$ is the radius of the wedge $W_1$, and $\alpha$ is its angle. Denote by $\chi_1(t)$ a smooth cut-off function that is identically 1 for $0 \leq t \leq r_0/2$, and identically 0 for $t \geq 3r_0/4$, and by $\chi_2(t)$ a cut-off function that has these properties, but with $r_1$ instead of $r_0$. Then the forms $\alpha_j := \overline{\partial}(f_j(z_2)\chi_1(|z_1|)\chi_2(|z_2|))$ are $\overline{\partial}$-closed on $\Omega$. Let $g_j := \overline{\partial}^* N_1 \alpha_j$. After passing to a subsequence, we may assume that $\{g_j\}_{j=1}^\infty$ converges in $L^2(\Omega)$ if $\overline{\partial}^* N_1$ is compact. Set $h_j(z_1, z_2) := \chi_1(|z_1|)\chi_2(|z_2|)f_j(z_2) - g_j(z_1, z_2)$. Then the $h_j$ are holomorphic in $\Omega$, and $h_j = -g_j$ outside the support of $\chi_1(|z_1|)\chi_2(|z_2|)$ (and hence $\{h_j\}_{j=1}^\infty$ converges in $L^2$ of $\Omega$ minus that support). Using the mean value property of holomorphic functions in the $z_1$ variable gives that $\{h_j\}_{j=1}^\infty$ also converges in $L^2(D_{r_0/2} \times (W_0 \cap D_{r_3}))$. This convergence also holds for $\{g_j\}_{j=1}^\infty$, hence for $\{\chi_1(|z_1|)\chi_2(|z_2|)f_j(z_2)\}_{j=1}^\infty$. But this implies that $\{f_j\}_{j=1}^\infty$ converges in $W_0 \cap D_{r_3}$ for any $r_3 < r_1/2$, a contradiction. Therefore, $\overline{\partial}^* N_1$ is not compact.

The general case of Proposition 9 is obtained by replacing the products disc × wedge above by suitable biholomorphic images of such products: the crucial noncompactness of the restriction operator between the Bergman spaces is invariant under these biholomorphisms.

Whether there can be obstructions to compactness (say for $(0, 1)$-forms) more subtle than discs in the boundary has been settled only recently. Peter Matheos [M97] showed that there exist smoothly bounded complete Hartogs domains in $\mathbb{C}^2$ without discs in their boundaries, but whose $\overline{\partial}$-Neumann operators are nonetheless not compact. Such a domain $\Omega$ will be of the form $\Omega = \{(z, w) \in \mathbb{C}^2; \ z \in \Omega_1, |w| < e^{-\Phi(z)}\}$, where $\Omega_1$ is a domain in $\mathbb{C}$ and $\Phi$ is smooth and subharmonic in

$\Omega_1$. Assume that the boundary points of the form $(z, 0)$ are strictly pseudoconvex. The weakly pseudoconvex points are the points $\{(z,w) \in b\Omega; \; \triangle\Phi(z) = 0\}$. It is easy to see that $b\Omega$ will contain an analytic disc if and only if the set $W := \{z \in \Omega_1; \; \triangle\Phi(z) = 0\}$ has non-empty interior (to produce an analytic disc, consider a conjugate harmonic function to $\Phi$ on a (small) disc contained in the interior of $W$). The boundary will fail to satisfy property $(P_1)$ (equivalently: $B$-regularity) if and only if $W$ has non-empty *fine* interior, see [Si87], p. 310. Recall that the fine topology is the smallest topology that makes all subharmonic functions continuous; see, *e.g.*, [He69] for properties of this topology. It is strictly larger than the Euclidean topology, and there exist compact sets with empty Euclidean interior, but non-empty fine interior. Matheos [M97] constructed examples of compact sets $W$ with empty Euclidean interior and associated complete Hartogs domains whose weakly pseudoconvex boundary points project onto $W$ and whose $\overline{\partial}$-Neumann operator $N_1$ is not compact. (Such sets $W$ then necessarily have nonempty fine interior, by Sibony's result). Matheos' ideas combine with results from potential theory to give the following: *Every* compact set $W \subset \mathbb{C}$ with empty Euclidean interior but non-empty fine interior arises in this way.

**Theorem 10.** *Assume that $W$ is a compact subset of $\mathbb{C}$ with non-empty fine interior. Then there exists a smooth, bounded, pseudoconvex, complete Hartogs domain $\Omega$ in $\mathbb{C}^2$ whose weakly pseudoconvex boundary points project onto $W$ and whose $\overline{\partial}$-Neumann operator $N_1$ is not compact.*

**Remark 5.** In particular, when $W$ has empty Euclidean interior, the complete Hartogs domain will have no discs in the boundary. Note that it does not follow from Theorem 10 that if the boundary of a smooth Hartogs domain does not satisfy property $(P_1)$ (equivalently: the associated set $W$ in the plane has non-empty fine interior), then its $\overline{\partial}$-Neumann operator $N_1$ is not compact. Whether or not this implication is true (that is, whether or not for complete smooth Hartogs domains in $\mathbb{C}^2$, property $(P_1)$ and compactness of $N_1$ are actually equivalent) remains open.

**Remark 6.** Although the Hartogs domains in Theorem 10 have non-compact $\overline{\partial}$-Neumann operator, this operator is regular in Sobolev spaces; see [BS92], Theorem 1: any bounded smooth Hartogs domain in $\mathbb{C}^2$ with no disc in the boundary is "nowhere wormlike", in the terminology of [BS92], and so has globally regular $\overline{\partial}$-Neumann operator $N_1$.

Matheos first uses the machinery of Chang, Nagel, and Stein [CNS92] to show that a compactness estimate on $\Omega$ is equivalent to a compactness estimate on the boundary, in terms of $\overline{\partial}_b$ and $\overline{\partial}_b^*$. This equivalence is established for all smooth pseudoconvex domains in $\mathbb{C}^2$. Then, the key property of the set $W$ is that there exists a sequence of "Rayleigh functions" with supports shrinking to $W$. This last property turns out to hold for every compact set with non-empty fine interior ([F99], section 3). We now indicate how to combine this fact with Matheos' arguments to prove Theorem 10. We work on $\Omega$ directly; this results in some simplifications.

For an open set $\Omega$ in $\mathbb{R}^n$, let $\lambda(\Omega) = \inf\{\|\nabla f\|^2; \; f \in W_0^1(\Omega), \|f\| = 1\}$. The classical Rayleigh-Ritz formula says that $\lambda(\Omega)$ is the first eigenvalue of the Dirichlet problem for the (negative) Laplacian if $\Omega$ is regular with respect to the Dirichlet problem (see, e.g., [Cha84]). The Dirichlet problem has also been formulated for finely open sets in $\mathbb{R}^n$ (cf. [F72]). Among many properties, the monotonicity of $\lambda$ with respect to the domain remains true for finely open sets. That is, if $\Omega_2 \supset \Omega_1 \neq \emptyset$

are finely open sets in $\mathbb{R}^n$, then $\lambda(\Omega_2) \leq \lambda(\Omega_1) < \infty$. We refer the reader to [F99], section 3, for a discussion of these (and further) results concerning the Dirichlet Laplacian on finely open sets.

Without loss of generality, we assume that $W$ is contained in the open unit disc $D$. Let $\{D_j\}_{j=1}^{\infty}$ be the connected components of $D \setminus W$. Let $W_k = D \setminus \left(\cup_{j=1}^{k} \overline{D}_j\right)$. It follows from Theorem 10.14 in [He69] that $W_k \supset \text{int}_f(W) \neq \emptyset$, where $\text{int}_f$ denotes the interior in the fine topology. Therefore $\lambda(W_k) \leq \lambda(\text{int}_f(W)) < \infty$. It follows that there exists a (Rayleigh) sequence of functions $v_k \in C_0^{\infty}(W_k)$ such that $\|v_k\| = 1$ and $\|\nabla v_k\| \lesssim 1$. (Throughout the rest of this section, $A(k) \lesssim B(k)$ means that $A(k) \leq C \cdot B(k)$, where $C$ is some positive constant independent of $k$.)

Let $\varphi$ be a non-negative smooth function on $\mathbb{C}$ that vanishes to infinite order on $W$ and is strictly positive on $\mathbb{C} \setminus W$. Let $\psi_j = \varphi$ on $\overline{D}_j$ and $\psi_j = 0$ on $D \setminus \overline{D}_j$. Then $\varphi = \sum_{j=1}^{\infty} \psi_j$ on $\overline{D}$ and $\psi_j \in C^{\infty}(\overline{D})$. We now choose by induction a sequence $\{c_j\}_{j=1}^{\infty}$ in $(0, 1]$ and a strictly increasing sequence of positive integers $\{n_j\}_{j=1}^{\infty}$ such that the functions $\widetilde{\psi}_j = c_j \psi_j$ satisfy the following: (1) $n_j \int_{D_j} \widetilde{\psi}_j dA = 2\pi$; (2) $\|\widetilde{\psi}_{j+1}\|_{L^{\infty}} \leq \|\widetilde{\psi}_j\|_{L^{\infty}}$; (3) $n_j \|\widetilde{\psi}_{j+1}\|_{L^{\infty}} \leq 1$; (4) $n_{j+1}$ is divisible by $n_j$. First choose $n_1$ sufficiently large so that $c_1 = 2\pi/(n_1 \int_{D_1} \psi_1 dA) \leq 1$. Suppose that we have chosen $c_j$ and $n_j$. We can choose an integer $m \geq 2$ sufficiently large so that $c_{j+1} = 2\pi/(m \cdot n_j \int_{D_{j+1}} \psi_{j+1} dA) \leq \min\{1, \|\widetilde{\psi}_j\|_{L^{\infty}}/\|\psi_{j+1}\|_{L^{\infty}}, 1/(n_j\|\psi_{j+1}\|_{L^{\infty}})\}$. We then let $n_{j+1} = m \cdot n_j$.

We now let $\widetilde{\varphi} = \sum_{j=1}^{\infty} \widetilde{\psi}_j$. Then $\widetilde{\varphi} \in C^{\infty}(\overline{D})$. (The supports of $\widetilde{\psi}_j = c_j \psi_j$ are disjoint, except for points in $W$, where *all* derivatives vanish. Consequently the series actually converges in, for example, Sobolev norms.) Let $G = (1/2\pi) \log |z|$ and let $\Phi(z) = \int_D G(z - w) \widetilde{\varphi}(w) dA(w)$. Then $\Phi(z) \in C^{\infty}(\overline{D})$ and $\triangle \Phi = \varphi$. Note that $\Phi$ is strictly subharmonic near $\partial D$. Extend $\Phi$ to the disc $\widetilde{D}$ with center $0$ and radius $2$ such that $\Phi$ is smooth on $\widetilde{D}$, strictly subharmonic on $\widetilde{D} \setminus D$, and equals $-(1/2) \log(4 - |z|^2)$ near $\partial \widetilde{D}$. (Extend $\Phi$ smoothly to a function $\Phi_1$ compactly supported in $\widetilde{D}$. Then add a radial function that is identically $0$ for $|z| \leq 1$ and agrees with $-(1/2)\log(4-|z|^2)$ when $|z|$ is close to $2$. Such a function can be chosen to have its second (radial) derivative as big as we wish on a given compact subset of $\{1 < |z| < 2\}$, in particular on $\{\Delta \Phi_1 \leq 0\}$.) Then $\Omega = \{(z, w); \ z \in \widetilde{D}, \ |w| < \exp(-\Phi(z))\}$ is smooth. It remains to see that the $\overline{\partial}$-Neumann operator of $\Omega$ is not compact. For this, we shall construct a sequence of $(0, 1)$-forms that contradicts a compactness estimate, *i.e.*, we construct a sequence $\{u_k\}_{k=1}^{\infty}$ in $\text{Dom}(\overline{\partial}^*) \cap C^{\infty}(\overline{\Omega})$ such that the estimate

$$(8) \qquad \|u_k\|^2 \leq \varepsilon(\|\overline{\partial} u_k\|^2 + \|\overline{\partial}^* u_k\|^2) + C_{\varepsilon} \|u_k\|_{-1}^2,$$

uniformly in $k$, will fail for some $\varepsilon > 0$.

Let $\widetilde{\varphi}_k = \sum_{j=1}^{k} \widetilde{\psi}_k$ and let $\Phi_k = \int_D G(z-w) \widetilde{\varphi}_k(w) dA(w)$. Then $\triangle \Phi_k = \widetilde{\varphi}_k = 0$ on $W_k$. It follows from the way we construct $\widetilde{\psi}_j$ that $n_k \Phi_k$ has a harmonic conjugate $\Theta_k$ on $W_k$ whose values are determined up to an integer multiple of $2\pi$. Therefore $\exp(i\Theta_k)$ is single valued and smooth on $W_k$. Let $f_k = \sqrt{n_k} v_k(z) \exp(n_k \Phi + i\Theta_k) w^{n_k - 1}$ and let $u_k = f_k \cdot (\bar{w} d\bar{z} - 2\exp(-2\Phi)\Phi_z d\bar{w})$. Then $u_k \in \text{Dom}(\overline{\partial}^*) \cap C^{\infty}(\overline{\Omega})$ and $\|u_k\|^2 \approx \|f_k\|^2 \approx 1$. Moreover, $\{f_k\}$ and hence $\{u_k\}$ converge weakly to zero in $L^2(\Omega)$, and $L^2_{(0,1)}(\Omega)$, respectively (note that $f_j \perp f_k$ in $L^2(\Omega)$ when $j \neq k$) and so $\{u_k\}$ converges to zero in $W_{(0,1)}^{-1}(\Omega)$ in norm. (Alternatively, one easily checks that $\|u_k\|_{-1}^2 \approx \|f_k\|_{-1}^2 \lesssim (1/n_k^2)\|f_k\|^2 \approx 1/n_k^2$ .) Therefore, in order to see that estimate (8) is violated for some $\epsilon > 0$, it suffices to see that $\|\overline{\partial}u_k\|^2 + \|\overline{\partial}^* u_k\|^2$ is bounded independently of $k$. The Kohn-Morrey formula gives

$$\|\overline{\partial}u_k\|^2 + \|\overline{\partial}^* u_k\|^2 = \int_{b\Omega} \sum_{j,\ell=1}^2 \frac{\partial^2 \rho}{\partial z_j \partial \bar{z}_\ell}(u_k)_j \overline{(u_k)_\ell} d\sigma + \int_\Omega \sum_{j,\ell=1}^2 \left|\frac{\partial(u_k)_\ell}{\partial \bar{z}_j}\right|^2 dV,$$

where (temporarily) $z_1 = z$, $z_2 = w$. $\rho$ is a defining function normalized so that $|\nabla \rho| = 1$ on $b\Omega$. We may take this function to be $w\bar{w} - \exp(-2\Phi(z))$, suitably normalized. Then the right hand side is dominated by a positive constant (independent of $k$) times

$$\int_{\partial\Omega} \Phi_{z\bar{z}} |f_k|^2 d\sigma + \int_\Omega \left(|f_k|^2 + |f_{k\bar{z}}|^2\right) dV$$
$$\lesssim \int_{W_k} \left(n_k |v_k|^2 \Phi_{z\bar{z}} + |v_k|^2 + |v_{k\bar{z}}|^2 + |v_k (n_k \Phi + i\Theta_k)_{\bar{z}}|^2\right) dA.$$

The integral of the first term in the last expression is not greater than $n_k \|\widetilde{\varphi}\|_{L^\infty(W_k)} \times \int_{W_k} |v_k|^2 dA = n_k \|\widetilde{\psi}_{k+1}\|_{L^\infty(W_k)} \leq 1$. The contribution from the second and third terms is under control because of our choice of the sequence $\{v_k\}_{k=1}^\infty$. The integral of the last term is not greater than $\|(n_k \Phi + i\Theta_k)_{\bar{z}}\|_{L^\infty(W_k)}^2$. Since $(n_k \Phi + i\Theta_k)_{\bar{z}} = n_k (\Phi - \Phi_k)_{\bar{z}}$ for $z \in W_k$ (because $(n_k \Phi_k + i\Theta_k)_{\bar{z}} = 0$), and the modulus of the latter term is not greater than

$$\frac{n_k}{4\pi} \int_D \frac{|\widetilde{\varphi} - \widetilde{\varphi}_k|}{|z - w|} dA(w) \lesssim n_k \|\widetilde{\varphi} - \widetilde{\varphi}_k\|_{L^\infty(D)} = n_k \|\widetilde{\psi}_{k+1}\|_{L^\infty(W_k)} \leq 1,$$

it follows that $\|\overline{\partial}u_k\|^2 + \|\overline{\partial}^* u_k\|^2 \lesssim 1$. This concludes the proof of Theorem 10.

## 5. Domains related to convex domains

The sufficient conditions for compactness discussed in Section 3 are also necessary in the case of convex domains, and, more generally, in the case of domains that are locally convexifiable. Moreover, for these domains, these conditions are equivalent to the absence of (germs of) varieties from the boundary. Pseudoconvex Reinhardt domains are "almost" in this latter class, that is, they are locally convexifiable at the boundary points away from the coordinate hyperplanes (but may not be so at boundary points on the coordinate hyperplanes, see, *e.g.*, [FIK96]), and some of the results carry over to this class.

We say that a domain is locally convexifiable if for every boundary point there is a neighborhood, and a biholomorphic map defined on this neighborhood, that takes the intersection of the domain with the neighborhood onto a convex domain. Note that boundaries of convex domains are locally graphs of Lipschitz functions, so the same is true for the boundary of a locally convexifiable domain. The following theorem comes from [FS98].

**Theorem 11.** *Let $\Omega$ be a bounded pseudoconvex domain in $\mathbb{C}^n$ which is locally convexifiable, let $1 \leq q \leq n$. The following are equivalent:*
1. *The $\overline{\partial}$-Neumann operator $N_q$ is compact.*
2. *The boundary of $\Omega$ does not contain any analytic variety of dimension $\geq q$.*

3. *The boundary of $\Omega$ satisfies property $(P_q)$.*

Note that for a locally convexifiable domain, property $(P_q)$ is equivalent to $B_q$-regularity. Henkin and Iordan [HI97] had shown earlier that in the case of convex domains, if the boundary contains no analytic varieties of dimension $\geq 1$, then all $\overline{\partial}$-Neumann operators $N_q$, $1 \leq q \leq n$, are compact.

Actually, that (2) $\Rightarrow$ (3) holds for the general case of locally convexifiable domains is only implicit in [FS98]. We now indicate how to make this explicit. By Proposition 7, property $(P_q)$ for $b\Omega$ will follow if we show that for any $z \in b\Omega$, the point mass is the only $P_q$-measure on $b\Omega$ at $z$. This last property will follow if we can establish it locally, again by Proposition 7. So, let now $K$ be a compact subset of the boundary, small enough so that it is contained in some open set on which there is a biholomorphism that convexifies the boundary, say $w = f(z) = (f_1(z), \ldots, f_n(z))$. By shrinking $K$ if necessary, we may assume that there is a (small) domain $U$ such that $K \subset U$, $f(U \cap \Omega)$ is convex, $f$ is biholomorphic in a neighborhood of $\overline{U \cap \Omega}$, and there is no variety of dimension greater than or equal to $q$ in the boundary of $f(U \cap \Omega)$. By [FS98], Proposition 3.2, there exists, through each boundary point, a complex affine subspace $L$ of dimension $\leq q-1$ such that $L \cap \overline{f(U \cap \Omega)}$ is a peak set for $A(f(U \cap \Omega))$, the algebra of functions analytic in $f(U \cap \Omega)$ and continuous on $\overline{f(U \cap \Omega)}$. The inverse image under $f$ of $L$ gives a manifold of dimension $\leq q-1$ such that $f^{-1}(L) \cap \overline{U \cap \Omega}$ is a peak set for $A(U \cap \Omega)$. For every boundary point in $K$, there exists a peak set for $A(U \cap \Omega)$ of this form that contains the point. Let $g$ be the corresponding peak function. Note that every element in $A(U \cap \Omega)$ can be approximated uniformly on $\overline{U \cap \Omega}$ by functions analytic in a neighborhood of $\overline{U \cap \Omega}$, since this is true for the corresponding algebra on the *convex* domain $f(U \cap \Omega)$, and $f$ is biholomorphic in a neighborhood of $\overline{U \cap \Omega}$. Thus $|g| \in P_q(K)$. Consequently, for $z \in K$, all $P_q$-measures of $K$ at $z$ are supported on $K \cap f^{-1}(L)$. Let now $z_0 \in K$, and choose coordinates in $\mathbb{C}^n$ such that $f(z_0) = 0$ and $L \subset \{w_q = w_{q+1} = \ldots = 0\}$. Then for sufficiently large $M$, the function $-\sum_{j=1}^{q-1} |f_j(z)|^2 + M \sum_{j=q}^{n} |f_j(z)|^2$ is in $P_q(K)$. It equals zero at $z_0$, and is non-positive on $K \cap f^{-1}(L)$. If we integrate against a $P_q$-measure for $z_0$ over $K$, and take into account that the measure is supported on $K \cap f^{-1}(L)$, it follows that it is actually a point mass at $z_0$.

**Remark 7**. The proof of (1) $\Rightarrow$ (2) in [FS98] establishes the following: if $N_q$ is compact and $\Omega$ is convexifiable at (hence near) a boundary point $p$, then the boundary cannot contain $q$-dimensional varieties near $p$. The reduction to the (globally) convex case uses the localization in Lemma 2 in an essential way, together with the observation that compactness of the relevant $\overline{\partial}$-Neumann operators is invariant under the convexifying biholomorphism (Lemma 1, (4), and the observation that compactness of the canonical solution operator follows once there is some compact solution operator), see [FS98], Remark (1) in section 5.

**Remark 8**. On *convex* domains, more is true. As mentioned in Section 2, compactness of the commutators $[P_{q-1}, \bar{z}_j]$ implies compactness of $\overline{\partial}^* N_q$ restricted to forms with holomorphic coefficients (always), and this in turn implies compactness of $\overline{\partial}^* N_q$ on all of $L^2_{(0,q)}(\Omega)$ if $\Omega$ is convex ([FS98], Remark (2), Section 5). Thus in the convex case, properties (1)–(3) in Theorem 11 are furthermore equivalent to compactness of the commutators $[P_{q-1}, M]$ discussed in Proposition 4 as well as to compactness of $\overline{\partial}^* N_q$, restricted to forms with holomorphic coefficients. In addition, the absence of varieties from the boundary is equivalent with the absence of affine varieties. This is a simple manifestation of the general principle that on convex boundaries, questions of orders of contact of varieties are decided by orders of contact of affine varieties. (See [FS98], Section 2, and the references given there.)

As mentioned earlier, pseudoconvex Reinhardt domains are locally convexifiable at most of their boundary points, so it is not surprising that the above ideas also give results on this class of domains.

**Theorem 12.** *Let $\Omega$ be a bounded pseudoconvex Reinhardt domain in $\mathbb{C}^n$, $1 \leq q \leq n$. If the boundary of $\Omega$ does not contain an analytic variety of dimension $\geq q$, then it is $B_q$-regular, and consequently, the $\overline{\partial}$-Neumann operator $N_q$ on $(0, q)$-forms is compact.*

The fact that if the boundary of a Reinhardt domain contains no varieties of dimension $\geq 1$, then $N_q$, $1 \leq q \leq n$, is compact is in [HI97], Corollary 1. However, there is an inaccuracy in the proof: it is asserted that on a bounded pseudoconvex Reinhardt domain $\Omega$, every boundary point is a peak point for the algebra $A(\Omega)$, the algebra of continuous functions on $\overline{\Omega}$ which are holomorphic on $\Omega$. Pseudoconvex Reinhardt domains in $\mathbb{C}^2$ lying over the punctured disc show that this is not the case in general: the boundary point corresponding to the puncture is not a peak point. Likewise, such a domain will not satisfy property $(B)$ from [HI97].

**Remark 9**.The example of the domain $\Omega = \{(z_1, z_2) \in \mathbb{C}^2;\ |z_1|^2 + |z_2|^2 < 1, 0 < |z_1| < 1\}$ discussed in Section 4 above shows that compactness of the $\overline{\partial}$-Neumann operator $N_q$ on a Reinhardt domain does in general *not* imply the absence of varieties from the boundary. It does imply, in view of Remark 7, that *the boundary intersected with the complement of the coordinate hyperplanes cannot contain $q$-dimensional varieties if $N_q$ is compact.* (The implication that if $N_1$ is compact, then the boundary (away from the coordinate hyperplanes) cannot contain analytic discs is contained in work of Salinas concerning Toeplitz $C^*$-algebras on Bergman spaces. For this work, and further references, see the two surveys [Sa91, Sa95].)

We now prove Theorem 12. The proof is by induction on the dimension $n$, as follows. When $n = 1$, $\Omega$ is either a disk, a punctured disk, or an annulus. In all three cases, the boundary is $B$-regular. We assume now that the theorem holds for all dimensions $\leq n - 1$, and show that it then holds in dimension $n$. Since the domain $\Omega$ is pseudoconvex, the boundary is locally convexifiable near each boundary point $(z_1, \ldots, z_n)$ where none of the $z_j$ $(1 \leq j \leq n)$ are 0. As in the proof of Theorem 11 (the implication $(2) \Rightarrow (3)$), it follows that for every such boundary point $z$, there is $r > 0$ such that the compact set $b\Omega \cap \overline{B(z, r)}$ is $B_q$-regular. Consequently, the part of the boundary where all coordinates are non-zero (this is not a compact set) is a countable union of compact sets, all of which are $B_q$-regular.

The rest of the boundary is a finite union of compact sets of the form $b\Omega \cap \{z_{j_1} = \ldots = z_{j_k} = 0, z_{j_{k+1}} \neq 0, \ldots, z_{j_n} \neq 0\}$, where $\{j_1, j_2, \ldots, j_n\}$ is a permutation of $\{1, 2, \ldots, n\}$ and $k$ ranges from 1 to $n - 1$. (The case $k = n$ results in only one set, $b\Omega \cap \{0\}$, which is either empty or $\{0\}$. But the set $\{0\}$ is $B_q$-regular.) So fix $(j_1, \ldots, j_k)$. Assume first that $n - k \geq q$. We will show that $b\Omega \cap \{z_{j_1} = \ldots = z_{j_k} = 0, z_{j_{k+1}} \neq 0, \ldots, z_{j_n} \neq 0\}$ (which we may assume nonempty) is the boundary, in $\{z_{j_1} = \ldots = z_{j_k} = 0\} \equiv \mathbb{C}^{n-k}$, of the bounded Reinhardt domain $\Omega \cap \{z_{j_1} = \ldots = z_{j_k} = 0\}$. Since $\Omega$ is pseudoconvex, its logarithmic image, $\log |\Omega| := \{(\log |z_1|, \ldots, \log |z_n|);\ z \in \Omega, z_j \neq 0, 1 \leq$ $j \leq n\}$, is convex. Since there is a boundary point with a zero coordinate, it is unbounded. We use the tools and notions developed in Section 8 of [Ro70]. Let $\{z_\nu\}_{\nu=1}^\infty$ be a sequence of points in $\Omega$ that converges to a boundary point $(0,\ldots,0,z_{k+1},\ldots,z_n)$, with $z_{k+1} \neq 0,\ldots,z_n \neq 0$. (Without loss of generality, we assume that $(j_1,j_2,\ldots,j_n) = (1,2,\ldots,n)$.) It follows from the convexity of the unbounded set $\log|\Omega|$ (specifically, see [Ro70], Theorem 8.2 and Theorem 8.3) that there is a vector $r = (r_1,\ldots,r_k,0,\ldots,0) \in \mathbb{R}^n$ such that $\zeta+r$ also belongs to $\log|\Omega|$ for each $\zeta \in \log|\Omega|$ ($r$ is in the recession cone of $\log|\Omega|$ in the terminology of [Ro70]). This shows that the projection of $\Omega$ onto $\{z_1 = \ldots = z_k = 0\}$ consists entirely of points in $\overline{\Omega}$. Since the boundary of $\Omega$ contains no $q$-dimensional varieties, and since $n - k \geq q$, at least one of these points has to be an interior point of $\Omega$. Pseudoconvexity now implies that $\Omega$ is complete with respect to the variables $(z_1,\ldots,z_k)$. Therefore, the projection of $\Omega$ onto $\{z_1 = \ldots = z_k = 0\}$ equals the intersection of $\Omega$ with this this coordinate "plane". Consequently, $b\Omega \cap \{z_1 = \ldots = z_k = 0\}$ is the same as the boundary of $\Omega \cap \{z_1 = \ldots = z_k = 0\}$. This is the boundary of a pseudoconvex bounded Reinhardt domain in $\mathbb{C}^{n-k}$. It contains no $q$-dimensional varieties. By the induction assumption, it is $B_q$-regular in $\mathbb{C}^{n-k}$. By considering functions of the form $\lambda(z_{k+1},\ldots,z_n) + M(|z_1|^2 + \ldots + |z_k|^2)$, one sees that $b\Omega \cap \{z_1 = \ldots = z_k = 0\}$ also is $B_q$-regular as a subset of $\mathbb{C}^n$.

It remains to consider the case $n - k < q$. In this case, any compact subset of $\{z_1 = \ldots = z_k = 0\}$ is $B_q$-regular as a subset of $\mathbb{C}^n$ (consider functions in $P_q(\mathbb{C}^n)$ of the form $M(|z_1|^2 + \ldots + |z_k|^2)$).

We now know that $b\Omega$ is a countable union of compact, $B_q$-regular subsets of $\mathbb{C}^n$. Therefore, by Proposition 8, $b\Omega$ itself is $B_q$-regular as well.

**Remark 10.** Recall the second example from Section 4 above: $\Omega$ is the unit ball in $\mathbb{C}^2$ minus the variety $\{z_1 = 0\}$, and $N_1$ is compact on $\Omega$. This case is not covered by Theorem 12 (the boundary contains a disc). It is instructive to compare the two situations. In the above proof of Theorem 12, we have used the fact that $b\Omega$ contains no $q$-dimensional variety twice. Namely, we used it first to conclude that the projection of $\Omega$ onto $\{z_1 = \cdots = z_k = 0\}$ consists of interior points of $\Omega$. This use of the absence of varieties is merely convenient: one could "fill in" these points as in the example. However, in order to make the induction run, we need to know that the *boundary* of this projection is free of $q$-dimensional varieties, and this *is* crucial for the argument, unless $q = n - 1$. When $q = n - 1$, the latter property comes for free. From these ideas, there results the following characterization of compactness of $N_{n-1}$ on a bounded pseudoconvex Reinhardt domain $\Omega$ in $\mathbb{C}^n$: $N_{n-1}$ *is compact if and only if* $b\Omega \setminus \{z_1 \cdots z_n = 0\}$ *contains no* $(n-1)$-*dimensional variety.* For complete Reinhardt domains, this characterization was obtained in [SSU89].

## ACKNOWLEDGMENTS

Part of this work was done while the first author visited Texas A & M University in the summer of 1999. He thanks the Department of Mathematics and Harold and Heidi Boas for their hospitality.

16    SIQI FU AND EMIL J. STRAUBE

# References

[Be82]   Steven R. Bell, *The Bergman kernel function and proper holomorphic mappings*, Transactions of the American Mathematical Society **270** (1982), no. 2, 585–591.

[Bo88]   Harold P. Boas, *Small sets of infinite type are benign for the $\overline{\partial}$-Neumann problem*, Proceedings of the American Mathematical Society **103** (1988), no. 2, 569–578.

[BS92]   Harold P. Boas and Emil J. Straube, *The Bergman projection on Hartogs domains in $\mathbb{C}^2$*, Transactions of the American Mathematical Society **331** (1992), no. 2, 529–540.

[BS99]   \_\_\_\_\_\_, *Global regularity of the $\overline{\partial}$-Neumann problem: a survey of the $L^2$-Sobolev theory*, Several Complex Variables (M. Schneider and Y.-T. Siu, eds.), MSRI Publications, vol. 37, to appear.

[Ca81]   David Catlin, *Necessary conditions for subellipticity and hypoellipticity for the $\overline{\partial}$-Neumann problem on pseudoconvex domains*, Recent Developments in Several Complex Variables (John E. Fornæss, ed.), Annals of Mathematics Studies, no. 100, Princeton University Press, 1981, pp. 93–100.

[Ca83]   \_\_\_\_\_\_, *Necessary conditions for subellipticity of the $\overline{\partial}$-Neumann problem*, Annals of Mathematics (2) **117** (1983), no. 1, 147–171.

[Ca84a]  \_\_\_\_\_\_, *Global regularity of the $\overline{\partial}$-Neumann problem*, Complex Analysis of Several Variables (Yum-Tong Siu, ed.), Proceedings of Symposia in Pure Mathematics, no. 41, American Mathematical Society, 1984, pp. 39–49.

[Ca84b]  \_\_\_\_\_\_, *Boundary invariants of pseudoconvex domains*, Annals of Mathematics (2) **120** (1984), no. 3, 529–586.

[Ca87]   \_\_\_\_\_\_, *Subelliptic estimates for the $\overline{\partial}$-Neumann problem on pseudoconvex domains*, Annals of Mathematics (2) **126** (1987), no. 1, 131–191.

[CD97]   David Catlin and John D'Angelo, *Positivity conditions for bihomogeneous polynomials*, Math. Res. Lett. **4** (1997), 555-567.

[CNS92]  D.-C. Chang, A. Nagel, and E. M. Stein, *Estimates for the $\overline{\partial}$-Neumann problem in pseudoconvex domains of finite type in $\mathbb{C}^2$*, Acta Mathematica **169** (1992), 153–228.

[Cha84]  Isaac Chavel, *Eigenvalues in Riemannian geometry*, Pure and Applied Mathematics vol.115, Academic Press, 1984

[CS99]   So-Chin Chen and Mei-Chi Shaw, *Partial differential equations in several complex variables*, book to appear.

[Ch98]   Michael Christ, *Private communication*, 1998.

[Da82]   John D'Angelo, *Real hypersurfaces, orders of contact, and applications*, Annals of Mathematics (2) **115** (1982), no. 3, 615–637.

[DP81]   Klas Diederich and Peter Pflug, *Necessary conditions for hypoellipticity of the $\overline{\partial}$-problem*, Recent Developments in Several Complex Variables (John E. Fornæss, ed.), Annals of Mathematics Studies, no. 100, Princeton University Press, 1981, pp. 151–154.

[FK72]   G. B. Folland and J. J. Kohn, *The Neumann problem for the Cauchy-Riemann complex*, Annals of Mathematics Studies, no. 75, Princeton University Press, 1972.

[FIK96]  Siqi Fu, Alexander Isaev, Steven Krantz, *Finite type conditions on Reinhardt domains*, Complex Variables: Theory and Application **31** (1996), 357-363.

[FS98]   Siqi Fu and Emil J. Straube, *Compactness of the $\overline{\partial}$-Neumann problem on convex domains*, Jour. Functional Anal. **159**(1998), 629-641.

[F72]    Bent Fuglede, *Finely Harmonic Functions*, Lecture Notes in Math. vol. 289, Springer, 1972.

[F99]    \_\_\_\_\_\_, *The Dirichlet Laplacian on finely open sets*, Potential Anal. **10** (1999), 91–101.

[G78]    Theodore W. Gamelin, *Uniform algebras and Jensen measures*, London Mathematical Society Lecture Note Series, no. 32, Cambridge University Press, 1978.

[He69]   L. L. Helms, *Introduction to potential theory*, Wiley-Interscience, 1969.

[HI97]   G. M. Henkin and A. Iordan, *Compactness of the Neumann operator for hyperconvex domains with non-smooth B regular boundary*, Math. Ann. **307** (1997), 151-168.

[H65]    Lars Hörmander, *$L^2$ estimates and existence theorems for the $\overline{\partial}$ operator*, Acta Mathematica **113** (1965), 89–152.

[KN65]   J. J. Kohn and L. Nirenberg, *Non-coercive boundary value problems*, Communications on Pure and Applied Mathematics **18** (1965), 443–492.

[Kr88]   Steven G. Krantz, *Compactness of the $\overline{\partial}$-Neumann operator*, Proceedings of the American Mathematical Society **103** (1988), no. 4, 1136–1138.

[Kr92] \_\_\_\_\_\_, *Partial differential equations and complex analysis*, CRC Press, Boca Raton, FL, 1992.

[La54] Peter D. Lax, *Symmetrizable linear transformations*, Comm. Pure Appl. Math. **7** (1954), 633-647.

[Li85] E. Ligocka, *The regularity of the weighted Bergman projections*, pp. 197-203 in Seminar on deformations, J. Lawrynowicz, editors, Lecture Notes in Math. **1165**, Springer, 1985.

[M97] Peter Matheos, *A Hartogs domain with no analytic discs in the boundary for which the $\overline{\partial}$-Neumann problem is not compact*, Preprint, 1997. (To appear in Jour. of Geom. Anal.)

[Ra84] R. Michael Range, *The $\overline{\partial}$-Neumann operator on the unit ball in $\mathbb{C}^n$*, Math. Ann. **266** (1984), 449-456.

[Ro70] R. Tyrrell Rockafellar, *Convex analysis*, Princeton Mathematical Series, No. 28 Princeton University Press, Princeton, N.J. 1970.

[Sa91] Norberto Salinas, *Noncompactness of the $\overline{\partial}$-Neumann problem and Toeplitz $C^*$-algebras*, Several Complex Variables and Complex Geometry (Eric Bedford, John P. D'Angelo, Robert E. Greene, and Steven G. Krantz, eds.), vol. 3, Proceedings of Symposia in Pure Mathematics, no. 52, American Mathematical Society, 1991, Proceedings of the Thirty-seventh Annual Summer Research Institute held at the University of California, Santa Cruz, California, July 10–30, 1989, pp. 329–334.

[Sa95] \_\_\_\_\_\_, *Toeplitz $C^*$-algebras and several complex variables*, in Multivariable operator theory, Contemp. Math. **185**, Amer. Math. Soc., Providence RI, 1995.

[SSU89] N. Salinas, A. Sheu, and H. Upmeier, *Toeplitz operators on pseudoconvex domains and foliation $C^*$-algebras*, Ann. of Math. **130** (1989), 531–565.

[Si87] N. Sibony, *Une classe de domaines pseudoconvexes*, Duke Mathematical Journal **55** (1987), no. 2, 299–319.

[Si89] \_\_\_\_\_\_, *Some aspects of weakly pseudoconvex domains*, Several Complex Variables and Complex Geometry (Eric Bedford, John P. D'Angelo, Robert E. Greene, and Steven G. Krantz, eds.), vol. 1, Proceedings of Symposia in Pure Mathematics, no. 52, American Mathematical Society, 1991, Proceedings of the Thirty-seventh Annual Summer Research Institute held at the University of California, Santa Cruz, California, July 10–30, 1989, pp. 119–232.

[St97] Emil J. Straube, *Plurisubharmonic functions and subellipticity of the $\overline{\partial}$-Neumann problem on non-smooth domains*, Math. Res. Lett. **4** (1997), 459-467.

[V72] U. Venugopalkrishna, *Fredholm operators associated with strongly pseudoconvex domains in $\mathbb{C}^n$*, Jour. Funct. Anal. **9** (1972), 349-373.

DEPARTMENT OF MATHEMATICS, UNIVERSITY OF WYOMING, LARAMIE, WY 82071
*E-mail address*: `sfu@uwyo.edu`

DEPARTMENT OF MATHEMATICS, TEXAS A & M UNIVERSITY, COLLEGE STATION, TX 77843
*E-mail address*: `straube@math.tamu.edu`